\documentstyle[12pt,fleqn]{article}

\font\twelve=cmbx10 at 15pt
\font\ten=cmbx10 at 12pt

\renewcommand{\theequation}{\arabic{section}.\arabic{equation}}

\newcommand{\ra}{\rightarrow}
\newcommand{\bra}{\langle}
\newcommand{\ket}{\rangle}
\newcommand{\be}{\begin{equation}}
\newcommand{\ee}{\end{equation}}
\newcommand{\bea}{\begin{eqnarray}}
\newcommand{\eea}{\end{eqnarray}}
\newcommand{\eps}{\varepsilon}

\newcommand{\ffi}{\varphi}
\newcommand{\ep}{$\spadesuit$}

\newtheorem{lem}{Lemma}[section]
\newtheorem{prop}{Proposition}[section]

\newtheorem{cor}{Corollary}[section]
\renewcommand{\thefootnote}{\alph{footnote}}
\def\R{\hbox{$\mit I$\kern-.277em$\mit R$}}
\def\C{\hbox{$\mit I$\kern-.7em$\mit C$}}
\def\un{\hbox{$\mit I$\kern-.77em$\mit I$}}

\begin{document}

\begin{titlepage}

\begin{center}

\renewcommand{\thefootnote}{\fnsymbol{footnote}}

{\ten Centre de Physique Th\'eorique\footnote{
Unit\'e Propre de Recherche 7061} - CNRS - Luminy, Case 907}
{\ten F-13288 Marseille Cedex 9 - France }

\vspace{1 cm}

\setcounter{footnote}{0}
\renewcommand{\thefootnote}{\arabic{footnote}}

{\twelve AN ADIABATIC THEOREM \\
FOR SINGULARLY PERTURBED HAMILTONIANS\footnote{Partially supported by
Fonds National Suisse de la Recherche Scientifique, Grant
8220-037200}}

\vspace{0.3 cm}

{\bf Alain JOYE\footnote{Permanent address: Centre de Physique
Th\'eorique, CNRS Marseille,
Luminy Case 907, F-13288 Marseille Cedex 9, France and PHYMAT,
Universit\'e de Toulon et du Var, B.P.132, F-83957 La Garde Cedex,
France
}}

\vspace{3,5 cm}

{\bf Abstract}

\end{center}

The adiabatic approximation in quantum mechanics is considered in the
case where the self-adjoint hamiltonian $H_0(t)$, satisfying the usual
spectral gap assumption in this context, is perturbed by
a term of the form $\eps H_1(t)$. Here $\eps\ra 0$ is the
adiabaticity  parameter and $H_1(t)$ is a self-adjoint operator
defined on a  smaller domain than the domain of $H_0(t)$. Thus the
total hamiltonian
$H_0(t)+\eps H_1(t)$ does not necessarily satisfy the gap assumption,
$\forall \eps >0$. It is shown that an adiabatic theorem can be
proven in  this situation under reasonnable hypotheses. The problem
considered can  also be viewed as the study of a time-dependent system
coupled to a time-dependent perturbation, in the limit of large
coupling constant.

\vspace{3,5 cm}


\bigskip

\noindent October 1994

\noindent CPT-94/P.3085

\bigskip

\noindent anonymous ftp or gopher: cpt.univ-mrs.fr

\end{titlepage}

\section{Introduction}

Due to its central role in Quantum Mechanics, the time-dependent
Schr\"odinger equation has been the object
of numerous studies. Since exact solutions to that equation are
rather  scarce, several asymptotic regimes governed by a set of
suitable parameters  were considered.
Among these asymptotic conditions, the so-called adiabatic regime is
a widely used limit in physics. The adiabatic limit describes the
evolution in time of a system when the governing hamiltonian is a
slowly varying function of time. A typical example is the slow, with
respect to the time scale of the system, switching on and off of time
dependent exterior perturbations. This regime also  plays an important
role in the study of models involving slow and fast  variables, like
molecular systems.  Mathematically speaking, the adiabatic limit is a
singular limit, as is the semi-classical limit. It corresponds to the
limit $\eps\ra 0$ of the equation
\be\label{schr}
  i\eps\psi'(t)=H_0(t)\psi(t),\;\;\; \psi(0)=\psi_0
\ee
where the prime denotes the time derivative, $\psi(t)$ is a vector
valued function in a separable Hilbert space
${\cal H}$ and $H_0(t)$ is a smooth enough family of self adjoint
operators, bounded from below and defined on some time independent
dense domain
$D_0$. We denote by $U_{\eps}(t)$ the associated unitary evolution
operator such that $U_{\eps}(0)=\un$. Let us further assume
that there is a gap in the spectrum of $H_0(t)$ for any
$t\in[0,1]$ and denote by $P_0(t)$ the spectral projector
corresponding to the bounded part of the spectrum. Then, the adiabatic
theorem of quantum mechanics asserts that there exists a unitary
operator $V(t)$ satisfying the intertwining property
\be\label{iw}
  V(t)P_0(0)=P_0(t)V(t),\;\;\;\forall t\in [0,1]
\ee
which approximates $U_{\eps}(t)$ in the limit $\eps\ra 0$
\be\label{atop}
  \sup_{t\in [0,1]}\|U_{\eps}(t)-V(t)\|={\cal O}(\eps).
\ee
This means in particular that
a solution of (\ref{schr}) such that
\be
  \psi(0)=P_0(0)\psi(0)
\ee
will satisfy
\be\label{atq}
  \psi(t)=P_0(t)\psi(t)+{\cal O}(\eps)
\ee
for any $t\in[0,1]$. In other words, such a solution follows the
spectral  subspace $P_0(t){\cal H}$ up to negligible corrections in
the limit
$\eps\ra 0$.
This formulation is a generalization of the early works \cite{bf} and
\cite{ka1} and can be found in \cite{n1}, \cite{asy}.
Recently, several refinements and further generalisations of this
result were  proven on a rigourous basis, see \cite{n2}, \cite{jp1}
and the references therein. Most of these works deal with the
construction of higher and higher  order approximate solutions to
(\ref{schr}) in the limit $\eps\ra 0$, even up to exponential order,
under  hypotheses similar to those loosely given above.
In particular, the spectral gap hypothesis is crucial for these
results to hold .
Another generalisation of the adiabatic theorem consists in trying to
get rid of the gap hypothesis. Such a result was obtained by Avron,
Howland and Simon
\cite{ahs} who
considered the case of hamiltonians with dense pure point spectrum.
Under certain conditions on the mismatch of resonnances they
could prove similar results to (\ref{atq}) for two kinds of
hamiltonians, with $P_0(t)$ the corresponding, either finite or
infinite dimensional,  spectral projector.

\subsection{The Problem}

In that paper we consider a related natural generalisation of the
adiabatic  theorem which consists in adding to the hamiltonian a
perturbation of order $\eps$.
We simply replace the hamiltonian $H_0(t)$ satisfying a gap
hypothesis by the perturbed
hamiltonian $H_0(t)+\eps H_1(t)$, where the self adjoint perturbation
$H_1(t)$ is defined on some time independent dense domain $D_1$.
This corresponds to the first order correction steming from a formal
hamiltonian $H(t,\eps)=H_0(t)+\eps H_1(t)+\eps^2 H_2(t)+\cdots$, for
example. If $D_1\supseteq D_0$,
regular pertubation theory shows that the total hamiltonian satisfies
the gap assumption for $\eps $ small enough, so that the usual
adiabatic theorem applies. However, if $D_1\subset D_0$, the term
$\eps H_1(t)$ becomes very singular eventhough $\eps\ra 0$ in the
adiabatic limit. In particular, if the gap hypothesis holds for
$H_0(t)$, it doesn't hold in general for the (suitable extension of)
the operator $H_0(t)+\eps H_1(t)$, no matter how small
$\eps$ is. We are thus in another case where the driving hamiltonian
has no gap in its spectrum.
We consider the equation
\be\label{scpe}
  i\eps U_{\eps}'(t)=(H_0(t)+\eps H_1(t))U_{\eps}(t)
\ee
on the suitable domain for the hamiltonian in the limit $\eps \ra 0$,
under the  main assumption that $H_0(t)$ satisfies the gap hypothesis.
Our main result, see Proposition \ref{mare}, is
the construction of an approximate evolution operator $V(t)$
satisfying (\ref{iw}) and (\ref{atop}), provided some regularity
conditions are satisfied. This result shows that the adiabatic theorem
survives some singular  perturbations or it can be viewed as providing
another situation where the  adiabatic theorem holds although
the driving hamiltonian does not satisfy the gap assumption. If the
hamiltonians $H_0(t)$ and $H_1(t)$ are both time-independent, a very
complete discussion of this problem can be found in \cite{n3}.

It should also be noticed that the equation (\ref{scpe}) rewritten as
\be\label{laco}
  iU_{\eps}(t)=(H_1(t)+\frac{1}{\eps}H_0(t))
  U_{\eps}(t),\;\;\;\eps\ra 0
\ee
can be viewed as describing the evolution of a system driven by the
hamiltonian $H_1(t)$ coupled to the perturbation $H_0(t)$ in the
limit  of infinite coupling constant.
We give below an application of our result in this setting.

\subsection{Heuristics}

We want to give here a formal argument explaining why this result
should hold and what the approximate evolution $V(t)$ should be.
Let $U_1(t)$ be defined by
\be
  iU_1'(t)=H_1(t)U_1(t), \;\;\; U_1(0)=\un
\ee
and let us consider the corresponding interaction picture.
The operator $\tilde{U}(t)\equiv U_1^{-1}(t)U_{\eps}(t)$ then
satisfies
\be
  i\eps\tilde{U}'(t)=U_1^{-1}(t)H_0(t)U_1(t)\tilde{U}(t)
  \equiv\tilde{H}(t)\tilde{U}(t).
\ee
The new hamiltonian $\tilde{H}(t)$ is $\eps$-independent and
satisfies a gap hypothesis if $H_0(t)$ does with corresponding
spectral projector given by
\be
  \tilde{P}(t)=U_1^{-1}(t)P_0(t)U_1(t).
\ee
In consequence, we can apply the standard adiabatic theorem to
$\tilde{U}(t)$.  Thus there exists an approximate evolution up to
order $\eps$, $\tilde{V}(t)$, of $\tilde{U}(t)$ such that
\be
  \tilde{V}(t)\tilde{P}(0)=\tilde{P}(t)\tilde{V}(t), \;\;\;\forall
t\in[0,1].
\ee
Coming back to the actual evolution this implies that
\be
  U_{\eps}(t)=U_1(t)\tilde{U}(t)= U_1(t)\tilde{V}(t)+{\cal O}(\eps)
\ee
where the approximation $V(t)\equiv U_1(t)\tilde{V}(t)$ is such that
\bea
  V(t)P_0(0)&=&U_1(t)\tilde{V}(t)\tilde{P}(0)\\
  &=&U_1(t)\tilde{P}(t)\tilde{V}(t)\nonumber\\
  &=&P_0(t)U_1(t)\tilde{V}(t)\nonumber\\
  &\equiv& P_0(0)V(t).\nonumber
\eea
It is known \cite{n1}, \cite{asy} that the approximate evolution
$\tilde{V}(t)$ is given by the solution of
\be
  i\eps \tilde{V}'(t)=\left(\tilde{H}(t)+i\eps[\tilde{P}'(t),
\tilde{P}(t)]
  \right)\tilde{V}(t), \;\;\; \tilde{V}(0)=\un.
\ee
Now, we compute
\be
  [\tilde{P}'(t),\tilde{P}(t)]=U_1^{-1}(t)\left([P_0'(t),P_0(t)]
  +i\left[[H_1(t),P_0(t)],P_0(t)\right]\right)U_1(t)
\ee
so that, $V(t)$ should satisfy
\bea
  i\eps V'(t)&=&(\eps H_1(t)+H_0(t)+i\eps [P_0'(t),P_0(t)]\\
  &-&\eps\left[[H_1(t),P_0(t)],P_0(t)\right])V(t) \nonumber\\
  &=&(H_0(t)+\eps P_0(t)H_1(t)P_0(t)+(\un -P_0(t))H_1(t)
  (\un -P_0(t)) \nonumber\\
  &+&i\eps [P_0'(t),P_0(t)])V(t).\nonumber
\eea
This equation is to be compared with (\ref{defv}). Of course, the
main problem to turn this sketch into a proof is the question of
domains. We want to stress the fact that our goal is to give here
reasonnable hypotheses  under which the adiabatic theorem holds,
although it should be possible to obtain higher order approximations
in this case too, provided additional assumptions are made.

In the next section we provide a precise statement of our main result
in  Proposition \ref{mare} and give a proof of it. Further remarks on
$V(t)$ and a  simple application are given at the end of the section.
The appendix  contains the proofs of technical lemmas.

\noindent{\bf Acknowledgements:} It is a pleasure to thank P.Briet,
P.Duclos, G.Hagedorn and E.Mourre for enlightening discussions.

\section{Main Result}
\subsection{Hypotheses}
\setcounter{equation}{0}
We start by expressing here
the hypotheses we need in order to prove our main proposition
\ref{mare}. \\ {\bf Hypothesis D:}

Let $H_0(t)$ and $H_1(t)$ be two time-dependent self-adjoint
operators in a separable Hilbert space ${\cal H}$ which are densely
defined on their   respective
domains $D_0$ and $D_1$ for all $t\in [0,1]$. These domains are
assumed to be independent of $t$. We introduce the operator
\be
  H(t,\eps)=H_0(t)+\eps H_1(t),\;\; D_0\cap D_1\ra {\cal H}
\ee
and we
further assume that there exists a dense domain $D\subseteq D_0\cap
D_1$ on which $H(t,\eps)$ is essentially self-adjoint. Thus there
exists a dense domain $D'\supseteq D_0\cap D_1\supseteq D$ on which
the closure
$\overline{H(t,\eps)}$ of $H(t,\eps)$ is self-adjoint.
This domain $D'$ is also supposed to be independent of $t$.\\
{\bf Hypothesis R$_0$:}

The operator
$H_0(t)$ is strongly $C^{2}$ on $D_0$ and its spectrum $\sigma_0(t)$
is divided into two parts $\sigma^{a}_0(t)$ and $\sigma^{b}_0(t)$
separated by a finite gap for all $t\in [0,1]$ such that
$\sigma^{b}_0(t)$ is bounded. Let $P_0(t)$ be the spectral projector
associated with the part
$\sigma^{b}_0(t)$ constructed by means of the Riesz formula
\be\label{depr}
  P_0(t)=-\frac{1}{2\pi i}\oint_{\Gamma}R_0(t,\lambda)d\lambda
\ee
where
$R_0(t,\lambda)=(H_0(t)-\lambda)^{-1}$ and $\Gamma$ is a path in the
resolvent set of $H_0(t)$ which encircles $\sigma^{b}_0(t)$. Note that
both $R_0(t,\lambda)$, $\lambda\in\Gamma$ and $P_0(t)$ are strongly
$C^{2}$ on ${\cal H}$. Moreover, for any $t\in [0,1]$,
\bea
  \mbox{Ran}P_0(t)&\subseteq&D_1\label{rep}\\
  \mbox{Ran}{P_0}'(t)&\subseteq&D_1\label{repp}
\eea
{\bf Hypothesis R$_1$:}
The operator $H_1(t)$ is strongly $C^{1}$
on $D_1$ and  $H_1(0){P_0}'(t)$ is strongly $C^{0}$ on  ${\cal H}$.
Moreover, for all $(\lambda ,t)\in \Gamma\times [0,1]$,
\bea \label{tecas}
  \mbox{Ran}R_0(t,\lambda)H_1(t)P_0(t)&\subseteq& D_1\\
  \mbox{Ran}R_0(t,\lambda)P_0'(t)P_0(t)&\subseteq& D_1\label{tecas1}\\
  \sup_{(\lambda ,t)\in\Gamma\times [0,1]}\| H_1(0)R_0(t,\lambda)
  H_1(t)P_0(t)\|&<&b<\infty \label{tecas2}\\
  \sup_{(\lambda ,t)\in \Gamma\times [0,1]}\|
  H_1(0)R_0(t,\lambda)P_0'(t)P_0(t)\|&<&b<\infty.\label{tecas3}
\eea
{\bf Hypothesis $\overline{\mbox{R}}$:}

We finally require the self adjoint operator
$\overline{H(t,\eps)}$ to be strongly $C^{1}$ on $D'$ and bounded
from  below for all $t\in [0,1]$.

\newpage
\noindent
{\bf Remarks:}

It follows from hypothesis D that $\overline{H(t,\eps)}$ converges
strongly to $H_0(t)$ as $\eps\ra 0$ in the generalized sense so that
the spectrum of $\overline{H(t,\eps)}$ is asymptotically concentrated
on any neighbourhood of the spectrum of $H_0(t)$ (\cite{ka2}, p.475).

Hypothesis D is of course also satisfied in the trivial case
$D_1\supseteq D_0$.

To check assumption (\ref{repp}), it is enough to show that
\be
  \mbox{Ran}{P_0}'(t){P_0}(t)\subseteq D_1
\ee
due to the identity
\be\label{proj}
  P_0(t)P_0'(t)P_0(t)\equiv 0.
\ee

It is also actually enough in $\overline{\mbox{R}}$ to require the
existence for  all $t\in [0,1]$ of a real
number in the resolvent set of $\overline{H(t,\eps)}$ instead of
asking
$\overline{H(t,\eps)}$ to be bounded from below.

\subsection{Preliminaries}

As a consequence of hypothesis $\overline{\mbox{R}}$, we have
a well defined unitary evolution $U_{\eps}(t)$  which together with
its inverse, is strongly differentiable on $D'$, maps $D'$ into $D'$
and
$U_{\eps}(t)$ satisfies
the Schr\"odinger equation
\be\label{esp}
  i\eps\frac{\partial}{\partial t}U_{\eps}(t)=\overline{\left(H_0(t)+
\eps
 H_1(t)\right)}U_{\eps}(t),\;\; U_{\eps}(0)=\un
\ee
(see \cite{rs}).  Another direct consequence of our hypotheses is the
following  technical lemma, the proof of which is given in appendix.

\begin{lem} \label{tele}
Under conditions D, R$_0$ and R$_1$, we have\\
a)\hspace{.6cm}$H_1(t)P_0(t)$ is bounded and
\bea
  P_0(t)H_1(t)&=& (H_1(t)P_0(t))^*\\
  H_1(t)P_0(t)&=&(P_0(t)H_1(t))^*
\eea
b)\hspace{.6cm}$H_1(t)P_0(t)$ is strongly $C^1$
on ${\cal H}$\\
c)\hspace{.6cm}$P_0(t)H_1(t)$ is strongly $C^1$ on ${\cal H}$.
\end{lem}
{\bf Remark:} We actually don't need conditions (\ref{tecas}) to
(\ref{tecas3}) to prove this result.

We set
\be
  H_1^{a}(t)=Q_0(t)H_1(t)P_0(t)+P_0(t)H_1(t)Q_0(t)
\ee
where $Q_0(t)=(\un -P_0(t))$. We have the immediate
\begin{cor}
The operator $H_1^{a}(t)$ is self-adjoint, bounded and strongly
$C^1$ on $\cal H$.
\end{cor}
We also introduce
\be
  K_0(t)=i[P_0'(t),P_0(t)]
\ee
which is bounded, self-adjoint and strongly $C^1$ as well.
Hence the perturbed operator
\be
  \overline{H_0(t)+\eps H_1(t)}-\eps H_1^{a}(t)+\eps K_0(t)
\ee
is self-adjoint, bounded from below and strongly $C^1$on $D'$ (see
\cite{ka2}). Thus there exists a
unitary operator $V(t)$ which together with its inverse maps $D'$ to
$D'$,  is strongly
differentiable on $D'$ and $V(t)$ satisfies
\be\label{defv}
  i\eps V'(t)=\left(\overline{H_0(t)+\eps H_1(t)}-\eps H_1^{a}(t)+
  \eps K_0(t)\right)V(t),\;\; V(0)=\un .
\ee
Let us check that
\be\label{com0}
  \left[\overline{H_0(t)+\eps H_1(t)}-\eps H_1^{a}(t),P_0(t)\right]
=0.
\ee
Consider $\ffi\in D\subseteq D_0\cap D_1$. Since $P_0(t)\ffi\in D_0
\cap D_1$  as well, we can write
\bea
  & &\left[\overline{H_0(t)+\eps H_1(t)}-\eps H_1^{a}(t),P_0(t)\right]
\ffi=\\
  & &\left[H_0(t)+\eps H_1(t)-\eps H_1^{a}(t),P_0(t)\right]\ffi=
\nonumber \\
  & &\eps\left[ H_1(t)-H_1^{a}(t),P_0(t)\right]\ffi=\nonumber \\
  & &\eps \left(H_1(t)P_0(t)-Q_0(t)H_1(t)P_0(t)\right)\ffi-
\nonumber\\
  & &\eps\left(P_0(t)H_1(t)-P_0(t)H_1(t)Q_0(t)\right)\ffi\equiv 0
\nonumber
\eea
using $Q_0(t)=\un -P_0(t)$.
Hence, (\ref{com0}) holds as $D$ is dense. It follows from
classical results (see \cite{kr}) that
\begin{lem}\label{intertwlem} Let $V(t)$ be defined by (\ref{defv}).
Then
\be\label{intertw}
  V(t)P_0(0)=P_0(t)V(t)
\ee
for all $t\in [0,1]$.
\end{lem}
This lemma shows that $V(t)$ follows the decomposition of the Hilbert
space
${\cal H}$ into ${\cal H}=P_0(t){\cal H}\bigoplus Q_0(t){\cal H}$.
Our goal  is now to show that this evolution is an approximation of
the actual evolution as $\eps \ra 0$:
\begin{prop}[Adiabatic Theorem]\label{mare} Assume D, $\mbox{R}_0$,
$\mbox{R}_1$ and
$\overline{\mbox{R}}$ and let $U_{\eps}(t)$ be defined by
(\ref{esp}). The  unitary $V(t)$ defined by (\ref{defv}) enjoying the
intertwining  property  (\ref{intertw}) is such that
\be
  \sup_{t\in[0,1]}\|U_{\eps}(t)-V(t)\|={\cal O} (\eps).
\ee
\end{prop}

\subsection{Proof of Proposition \protect\ref{mare}}

In order to compare $U_{\eps}(t)$ and $V(t)$, we introduce
the unitary operator $A(t)=V^{-1}(t)U_{\eps}(t)$. It maps $D'$ to
$D'$  and it satisfies for any $\ffi\in D'$
\bea\label{eda}
  iA'(t)\ffi&=&i{V^{-1}}'(t)U_{\eps}(t)\ffi+iV^{-1}(t)U_{\eps}'(t)
\ffi\\
  &=&V^{-1}(t)\left( H_1^{a}(t)-K_0(t)\right)V(t)A(t)\ffi, \;\; A(0)=
\un .
  \nonumber
\eea
We have used the unitarity of $V^{-1}(t)$ and the property
$U_{\eps}(t)\ffi\in D'$ to derive (\ref{eda}).
We want to perform an integration by parts on the Volterra equation
corresponding to (\ref{eda}), as in \cite{asy} and \cite{ahs}.
Let $B(t)$ be a bounded,
strongly $C^1$ operator and ${\cal R}(B)(t)$ be defined by
\be
  {\cal R}(B)(t)=\frac{1}{2\pi i}\oint_{\Gamma} R_0(t,\lambda) B(t)
R_0(t,\lambda) d\lambda,
\ee
where the path $\Gamma$ is as in (\ref{depr}).
\begin{lem}\label{GPR}
  If hypothesis R$_{0}$ holds, the operator ${\cal R}(B)(t)$ is
bounded, strongly $C^1$ and maps
${\cal H}$ into $D_0$. Moreover,
\bea
  a)\;\;\; [{\cal R}(B)(t),H_0(t)]&=&-[P_0(t),B(t)] \\
  b)\hspace{1.75cm} {\cal R}(B)(t)&=&P_0(t){\cal R}(B)(t)Q_0(t)+
Q_0(t){\cal R}(B)(t) P_0(t)\\
  c)\hspace{1.7cm}{\cal R}^*(B)(t)&=&{\cal R}(B^*)(t)
\eea
\end{lem}
The first identity is proven in \cite{asy} and the second, the proof
of which is given in appendix, is mentionned
in \cite{jp2}. The proof of the last assertion is straightforward and
will be  omitted. Hypotheses (\ref{rep}) and (\ref{repp}) are of
course not required to  get this lemma.

We need to control the range of ${\cal R}(B)(t)$ when
$B(t)=H_1^a(t)-K_0(t)$.
\begin{lem}\label{ran}
If hypotheses D, R$_0$ and R$_1$ hold, then the bounded operator
${\cal R}(H_1^a-K_0)(t)$
maps ${\cal H}$ into $D_0\cap D_1$. Moreover
\be
  \sup_{t\in [0,1]}\| H_1(t){\cal R}(H_1^a-K_0)(t)\|< \infty.
\ee
\end{lem}
{\bf Proof:} Consider the first statement. It follows from the first
assertion of lemma \ref{GPR} that it is
sufficient to show that $\mbox{Ran}{\cal R}(H_1^a-K_0)(t)\subseteq
D_1$ and  it follows from the part b)
of it that it is actually enough to look at
\be
  \mbox{Ran}Q_0(t){\cal R}(H_1^a-K_0)(t)P_0(t)=\mbox{Ran}
{\cal R}((H_1^a-K_0)P_0)(t).
\ee
Consider now the operator
\be\label{expb}
   H_1(t)R_0(t,\lambda) (H_1^a(t)-K_0(t))P_0(t) R_0(t,\lambda),\;\;\;
  \lambda\in \Gamma.
\ee
This operator is well defined because, see (\ref{proj}),
\bea
  & &R_0(t,\lambda)(H_1^a(t)-K_0(t))P_0(t)=\\
  & &R_0(t,\lambda)(Q_0(t)H_1(t)P_0(t)-iP_0'(t)P_0(t))=\nonumber\\
  & &R_0(t,\lambda)H_1(t)P_0(t)-P_0(t)R_0(t,\lambda)H_1(t)P_0(t)
  -iR_0(t,\lambda)P_0'(t)P_0(t)\nonumber
\eea
maps ${\cal H}$ into $D_1$ by hypothesis R$_1$. Moreover,
(\ref{expb}) is  uniformly bounded in $(\lambda,t)\in \Gamma\times
[0,1]$. Indeed,
\bea
  & &H_1(t)R_0(t,\lambda)(H_1^a(t)-K_0(t)) P_0(t)=\\
  & &H_1(t)(H_1(0)+i)^{-1}(H_1(0)+i)
  R_0(t,\lambda)(H_1^a(t)-K_0(t))P_0(t)\nonumber
\eea
where  $H_1(t)(H_1(0)+i)^{-1}$ is bounded due to the closed graph
theorem and strongly $C^1$ and
\bea
  & &\|H_1(0)R_0(t,\lambda)\left( H_1^a(t)-K_0(t)\right)P_0(t)
  \|\leq\\
  & &\|H_1(0)R_0(t,\lambda)H_1(t)P_0(t)\|+
  \|H_1(0)P_0(t)R_0(t,\lambda)
  H_1(t)P_0(t)\|+\nonumber\\
  & &\|H_1(0)R_0(t,\lambda)P_0'(t)P_0(t)\|\leq\nonumber\\
  & & 2b+\sup_{(\lambda,t)\in \Gamma\times [0,1]}\|H_1(0)P_0(t)\|
  \|R_0(t,\lambda)\|\|H_1(t)P_0(t)\|<\infty\nonumber
\eea
(since $H_1(0)P_0(t)$ is $C^1$, see (\ref{ap5}).)
This implies that both strong integrals
\bea\label{int}
  & &\frac{1}{2\pi i}\oint_{\Gamma} R_0(t,\lambda) (H_1^a(t)-K_0(t))
  P_0(t)R_0(t,\lambda)\psi d\lambda\\
  \label{int2}& &\frac{1}{2\pi i}\oint_{\Gamma} H_1(t)R_0(t,\lambda)
(H_1^a(t)-K_0(t))
  P_0(t)R_0(t,\lambda)\psi d\lambda
\eea
exist for any $\psi\in {\cal H}$. Hence, since
$H_1(t)$ is closed, the vector (\ref{int}) belongs to $D_1$ and
(\ref{int2}) equals $H_1(t)$ applied on (\ref{int}).
The second statement follows from the fact that
(\ref{int2}) is uniformly bounded in $t\in [0,1]$.\ep\\

We can now construct a modified integration by parts formula which is
the main  tool to prove the adiabatic theorem:
\begin{lem}\label{intl}
For any $\psi \in {\cal H}$, we can write
\bea
  & &Q_0(0)V^{-1}(s)(H_1^{a}(s)-K_0(s))V(s)P_0(0)\psi=\\
  & &i\eps\frac{d}{ds}\left(Q_0(0)V^{-1}(s){\cal R}(H_1^{a}-K_0)(s)
V(s)P_0(0)
  \psi\right)-\nonumber\\
  & &i\eps Q_0(0)V^{-1}(s){\cal R}'(H_1^{a}-K_0)(s)V(s)P_0(0)\psi
-\nonumber\\ & &\eps Q_0(0)V^{-1}(s)\left[{\cal
R}(H_1^{a}-K_0)(s),H_1(s)\right]V(s)P_0(0)
\psi\nonumber\\
& &\equiv i\eps\frac{d}{ds}\left(Q_0(0)V^{-1}(s){\cal R}(H_1^{a}-K_0)
(s)V(s) P_0(0)\psi\right)+\eps Q_0(0)C(s)P_0(0)\psi \nonumber
\eea
where $C(s)$ is bounded and satisfies
\be
  \sup_{s\in[0,1] \atop \eps \geq 0}\left\| C(s)\right\| <\infty.
\ee
\end{lem}
{\bf Proof:} On the one hand, using the definition of $H_1^{a}(s)$
and  the property (\ref{proj}), we can write using
lemmas \ref{intertwlem} and \ref{GPR}
\bea\label{hand}
  & &Q_0(0)V^{-1}(s)(H_1^{a}(s)-K_0(s))V(s)P_0(0)=\\
  & &Q_0(0)V^{-1}(s)[(H_1^{a}(s)-K_0(s)),P_0(s)]V(s)P_0(0)=
\nonumber\\
  & &Q_0(0)V^{-1}(s)[{\cal R}(H_1^{a}-K_0)(s),H_0(s)]V(s)P_0(0).
\nonumber
\eea
On the other hand, for any $\ffi\in D\subset D'$,
\bea\label{derv}
  & &i\eps\frac{d}{ds}\left(V^{-1}(s){\cal R}(H_1^{a}-K_0)(s)V(s)
\ffi\right)=\\
  & &i\eps {V^{-1}}'(s){\cal R}(H_1^{a}-K_0)(s)V(s)\ffi+
  i\eps V^{-1}(s){\cal R}(H_1^{a}-K_0)(s)V(s)'\ffi+\nonumber\\
  & &i\eps V^{-1}(s){\cal R}'(H_1^{a}-K_0)(s)V(s)\ffi =\nonumber\\
  & &V^{-1}(s)\left[{\cal R}(H_1^{a}-K_0)(s),
  \overline{H_0(s)+\eps H_1(s)}-\eps H_1^{a}(s)+\eps K_0(s)\right]
V(s)\ffi+\nonumber \\
  & &i\eps V^{-1}(s){\cal R}'(H_1^{a}-K_0)(s)V(s)\ffi\nonumber.
\eea
Note that
${\cal R}(H_1^{a}-K_0)(s)$ maps ${\cal H}$ into $D_0\cap D_1$
(lemma \ref{ran}), so that differentiation of the operator $V^{-1}
(s)$ on  the left of (\ref{derv}) is justified.
Replacing $\ffi$ by
$P_0(0)\psi$ and using
the property
\be
  V(s)P_0(0)\psi=P_0(s)V(s)P_0(0)\psi\in D_0\cap D_1\subset D'
\ee
we can expand the commutator
\bea\label{excom}
  & &\left[{\cal R}(H_1^{a}-K_0)(s),\overline{H_0(s)+\eps H_1(s)}
  -\eps H_1^{a}(s)+\eps K_0(s)\right]V(s)P_0(0)\psi =\\
  & &\left[{\cal R}(H_1^{a}-K_0)(s),H_0(s)+\eps H_1(s)
  -\eps H_1^{a}(s)+\eps K_0(s)\right]V(s)P_0(0)\psi=\nonumber\\
  & &\left[{\cal R}(H_1^{a}-K_0)(s),H_0(s)\right]V(s)P_0(0)
  \psi+\nonumber\\
  & &\eps\left[{\cal R}(H_1^{a}-K_0)(s),H_1(s)\right]V(s)P_0(0)
  \psi - \nonumber\\
  & &\eps\left[{\cal R}(H_1^{a}-K_0)(s),H_1^{a}(s)- K_0(s)
  \right]V(s)P_0(0)\psi.\nonumber
\eea
By lemma \ref{GPR}, part b) again, it is easily checked that
\be
  Q_0(s)\left[{\cal R}(H_1^{a}-K_0)(s),H_1^{a}(s)- K_0(s)
  \right]P_0(s)\equiv 0
\ee
so that the formula of lemma \ref{intl} follows. The commutator
\be
  \left[{\cal R}(H_1^{a}-K_0)(s),H_1(s)\right]
\ee
is uniformly bounded in $s$, as seen from lemma \ref{ran} and the
identity
\be
  {\cal R}(H_1^a-K_0)(t)H_1(s)=\left(H_1(s){\cal R}(H_1^a-K_0)(s)
\right)^*
\ee
which is proven as part a) of lemma \ref{tele}. Indeed, it is readily
checked form lemmas \ref{GPR} c) and \ref{ran} that ${\cal
R}(H_1^a-K_0)(s)$ has the  required properties.
The uniform boundedness in $s$ and $\eps$ of $C(s)$ then follows from
the  strong continuity
of ${\cal R}'(H_1^{a}-K_0)(s)$ and
from the unitarity of $V(s)$.\ep\\
{\bf Remark:} It is essential to consider $P_0(0)\psi$ instead of
$\ffi\in D$  in (\ref{excom}) to expand the operator
$\overline{H_0+\eps H_1}$  because $V(s)\ffi$ does not necessarily
belong to $D_0\cap D_1$.

Let us consider the projection in $Q_0(0){\cal H}$ of the integral
equation corresponding to (\ref{eda}).
\bea
  & &Q_0(0)A(t)\ffi -Q_0(0)\ffi=\\
  & &-i\int_0^t ds
  Q_0(0)V^{-1}(s)(H_1^{a}(s)-K_0(s))V(s)P_0(0)A(s)\ffi .\nonumber
\eea
We can apply lemma \ref{intl} and make use of (\ref{eda}) to write
\bea
  & &Q_0(0)A(t)\ffi -Q_0(0)\ffi = \\
  & &\eps\int_0^t ds \frac{d}{ds}\left(Q_0(0)V^{-1}(s){\cal R}
(H_1^{a}-K_0)(s)
  V(s)P_0(0)\right)A(s)\ffi-\nonumber \\
  & &i\eps \int_0^t ds Q_0(0)C(s)P_0(0)A(s)\ffi = \nonumber \\
  & & \left. \eps Q_0(0)V^{-1}(s){\cal R}(H_1^{a}-K_0)(s)V(s)
  P_0(0)A(s)\ffi\right|_0^t -\nonumber\\
  & &\eps \int_0^t ds Q_0(0)V^{-1}(s){\cal R}(H_1^{a}-K_0)(s)
  V(s)P_0(0)A'(s)\ffi-\nonumber\\
  & &i\eps \int_0^t ds Q_0(0)C(s)P_0(0)A(s)\ffi = \nonumber \\
  & &\left. \eps Q_0(0)V^{-1}(s){\cal R}(H_1^{a}-K_0)(s)V(s)
  P_0(0)A(s)\ffi\right|_0^t +\nonumber\\
  & &i\eps \int_0^t ds Q_0(0)V^{-1}(s){\cal R}(H_1^{a}-K_0)(s)
  (H_1^{a}(s)-K_0(s))V(s)Q_0(0)A(s)\ffi-\nonumber\\
  & &i\eps \int_0^t ds Q_0(0)C(s)P_0(0)A(s)\ffi.\nonumber
\eea
Since the only $\eps$-dependent operators $A(s)$ and $V(s)$ are
unitary and  the above operators are all uniformly bounded in $s\in
[0,1]$, we have arrived to
\be\label{esoa}
  \sup_{t\in[0,1]}\|Q_0(0)(A(t) -\un)\|={\cal O} (\eps).
\ee
Similarly,
\bea
  & &P_0(0)(A(t) - \un )P_0(0)\ffi=\\
  & &-i\int_0^t ds
  P_0(0)V^{-1}(s)(H_1^{a}(s)-K_0(s))V(s)Q_0(0)(A(s)-\un )P_0(0)
\ffi.\nonumber
\eea
 Taking (\ref{esoa}) into account we get
\be
  \sup_{t\in[0,1]}\|P_0(0)(A(t) -\un)P_0(0)\|={\cal O} (\eps),
\ee
hence
\be
  \sup_{t\in[0,1]}\|(A(t) -\un)P_0(0)\|={\cal O} (\eps).
\ee
Finally, since the projector
$P_0(0)$ is self-adjoint and $A(t)$ is unitary we can write
\bea
  \|P_0(0)(A(t) -\un)\|&=&\|(P_0(0)(A(t)-\un))^*\|\\
  &=&\|(A^{-1}(t)-\un)P_0(0)\|\nonumber\\
  &=&\|(\un -A(t))P_0(0)\|\nonumber\\
  &=&{\cal O} (\eps)\nonumber
\eea
uniformly in $t\in [0,1]$. This last equation and (\ref{esoa}) imply
$A(t)=\un + {\cal O}(\eps)$, uniformly in $t\in[0,1]$.
The definition of $A(t)$ eventually yields
\be
  \|(A(t) -\un)\|=\|V^{-1}(t)(U_{\eps}(t)-V(t))\|=\|(U_{\eps}(t)-V(t))
\|.
\ee
\hfill\ep

\subsection{Factorization of $V(t)$}

We can now decompose the approximate evolution $V(t)$ as in the usual
adiabatic theorem. Let $W(t)$ be defined by
\be
  iW'(t)=K_0(t)W(t)\; , \;\; W(0)=\un
\ee
and consider the unitary operator $\Phi(t)=W^{-1}(t)V(t)$. Since
$W(t)P_0(0)=P_0(t)W(t)$ for all $t$, we
have
\be\label{comm}
  [\Phi(t),P_0(0)]\equiv 0 ,\;\; \forall t\in [0,1].
\ee
The $\eps$-independent operator $W(t)$ has a geometrical meaning and
describes a parallel transport
of $P_0(t){\cal H}$ and $\Phi(t)$ is the analog of a dynamical phase
which  is singular in the limit $\eps\ra 0$. When restricted to
$P_0(0){\cal H}$, the operator $\Phi(t)$ satisfies a simple linear
differential equation with bounded generator.
\begin{lem}
Under hypotheses D, R$_0$, R$_1$ and $\overline{R}$, we have
\be
  V(t)P_0(0)=W(t)\Phi(t) P_0(0)=P_0(t)W(t)\Phi(t)P_0(0)
\ee
where $\Phi(t)P_0(0)$ satisfy
\be
  i\eps \Phi'(t)P_0(0)=W^{-1}(t)(H_0(t)P_0(t)
  +\eps P_0(t)H_1(t)P_0(t))W(t)\Phi(t)P_0(0).
\ee
\end{lem}
Accordingly,
If $P_0(0){\cal H}$ is one dimensional and $\psi =P_0(0)\psi$,
\be
  V(t)\psi=\exp\left\{-i\int_0^t(e_0(s)/\eps+e_1(s))ds\right\}W(t)
\psi,
\ee
where $e_0(t)$ is the associated eigenvalue of $H_0(t)$ and $e_1(t)$
is obtained by first order perturbation theory
\be
  e_1(t)=\mbox{Trace}P_0(t)H_1(t)P_0(t).
\ee
{\bf Proof:} Using (\ref{comm}), the intertwining property, and
hypothesis  R$_0$, we have
\bea
  i\eps\Phi'(t)P_0(0)&=& i\eps P_0(0)(W^{-1}(t)V(t))'P_0(0)\\
  &=&W^{-1}(t)P_0(t)(\overline{H_0(t)+\eps H_1(t)}-\eps H^a(t))
P_0(t)V(t)
  P_0(0)\nonumber\\
  &=&W^{-1}(t)P_0(t)(H_0(t)+\eps H_1(t))P_0(t)W(t)W^{-1}(t)V(t)
P_0(0).
  \nonumber
\eea \hfill \ep

\subsection{Example}

We present here a very simple application of proposition \ref{mare}.
Let ${\cal H}$ be the Hilbert space $L^2(\R ^n)$ and
\bea
  H_1(t)&=&-\frac{1}{2}\Delta+w(x,t)\;\;\;\mbox{on}\;\;D_1\\
  H_0(t)&=&\beta(t)|\ffi(t)\ket\bra\ffi(t)|\;\;\;
  \mbox{on}\;\; {\cal H}\nonumber
\eea
Here $w(x,t)$ is a smooth real valued function with compact support
in
$\R ^n\times [0,1]$
and $D_1$ is the domain of the self-adjoint extention of
$-\frac{1}{2}\Delta+w(x,t)$ on $C_0^{\infty}(\R ^n)$.
The unit vector $\ffi(t)\in{\cal H}$ for all $t\in [0,1]$ and
$\beta(t)$ is a real valued function of $[0,1]$.
We consider the equation
\be
  i\psi'(t)=\left(-\frac{1}{2}\Delta +w(x,t) +\frac{1}{\eps}
  \beta(t)|\ffi(t)\ket\bra \ffi(t)|\right)\psi(t)
\ee
on the domain $D'=D_1$
in the limit $\eps\ra 0$. This equation describes the dynamics of a
particule in a potential strongly coupled to a rank one operator, in
the  spirit of equation (\ref{laco}). We further assume that $\ffi(t)$
and
$\beta(t)$ are $C^2$ and that
\be
  \beta(t)\geq g>0\;,\;\;\forall t\in[0,1].
\ee
Thus the spectrum of $H_0(t)$,
$\sigma_0(t)=\{0,\beta(t)\}$, is divided into two disjoint parts.
Note that  the spectrum of $H_1(t)$ consists in the positive real axis
plus some negative eigenvalues, depending on $w(x,t)$. Since $H_0(t)$
is rank one, the spectrum of $H_0(t)+\eps H_1(t)$ is of the same
nature, so that $\beta(t)$ is not isolated, $\forall \eps >0$. Under
the additionnal hypotheses
\bea
  & &\Delta\ffi'(t)\in C^1[0,1] \\
  & &\sup_{t\in [0,1]}\|\Delta\Delta\ffi(t)\| <\infty\\
  & &\sup_{t\in [0,1]}\|\Delta w(\cdot ,t)\ffi(t)\| <\infty
\eea
it is readily checked that our hypotheses D, R$_0$, R$_1$ and
$\overline{\mbox{R}}$ are satisfied with
\be
  R_0(t,\lambda)=\frac{\beta(t)}{(\beta(t)-\lambda)\lambda}
  |\ffi(t)\ket\bra\ffi(t)|-\frac{1}{\lambda}\un
\ee
and
\be
  P_0(t)=|\ffi(t)\ket\bra\ffi(t)|.
\ee
By virtue of Proposition \ref{mare} and the above computations,
if $\psi(0)=\ffi(0)$, then
\bea
  \psi(t)&=&V(t)\ffi(0)+{\cal O}(\eps)\\
  &=&\exp\{-i\lambda(t,\eps)\}\ffi(t)+{\cal O}(\eps)\;, \;\;\forall
  t\in[0,1].\nonumber
\eea
Here
\bea
  & &\lambda(t,\eps)=\\
  & &\int_0^t(\beta(s)/\eps -\bra\ffi(s)|\Delta \ffi(s)\ket /2 +\bra
  \ffi(s)|w(\cdot ,s)\ffi(s)\ket - i\bra \ffi(s)|\ffi'(s)\ket)ds
\nonumber
\eea
is a real valued function and the last term in the integrant comes
from the parallel transport operator $W(t)$.

\appendix
\renewcommand{\theequation}{\Alph{section}.\arabic{equation}}

\section{Proof of Lemma \protect\ref{tele}:}
\setcounter{equation}{0}

a) We first show the boundedness of these operators. Since
$H_1(t)$ is self-adjoint, it is closed, and hypothesis (\ref{rep}) on
$P_0(t)$ implie that $H_1(t)P_0(t)$ is well defined on ${\cal H}$ and
is closed as well. Hence, by the closed graph theorem, it is  bounded.
Since all operators are densely defined, we have the general relation
\be\label{grel}
 (H_1(t)P_0(t))^* \supset P_0^*(t)H_1^*(t)
\ee
where $P_0^*(t)H_1^*(t)= P_0(t)H_1(t)$ and $H_1(t)P_0(t)$ is bounded.
As $D_1$ is dense in ${\cal H}$ we deduce from the extension
principle that
\be
  P_0(t)H_1(t)= (H_1(t)P_0(t))^*
\ee
is bounded.
Finally,
\be\label{adjb}
 (P_0(t)H_1(t))^* = H_1^*(t)P_0^{*}(t) =H_1(t)P_0(t)
\ee
holds since $P_0(t)$ is bounded, \cite{ka2} p.168. \\
b) Because of assumption (\ref{repp}) the operator
$H_1(0)P_0'(t)$ is bounded (see above) and it is strongly $C^0$ on
${\cal H}$ by hypothesis.
Hence the vector $H_1(0)P_0'(t)\ffi$, where $\ffi\in {\cal H}$, is
integrable and since $H_1(0)$ is closed we can write
\be
  \int_0^tH_1(0)P_0'(s)\ffi ds=H_1(0)\int_0^tP_0'(s)\ffi ds
  =H_1(0)(P_0(t)-P_0(0))\ffi .
\ee
Thus the bounded operator $H_1(0)P_0(t)$ is strongly $C^1$ and
\be\label{ap5}
  (H_1(0)P_0(t))'=H_1(0)P_0'(t).
\ee
Consider
\be
  H_1(t)P_0(t)=H_1(t)(H_1(0)+i)^{-1}(H_1(0)+i)P_0(t)
\ee
where $H_1(t)(H_1(0)+i)^{-1}$ is also bounded and strongly
$C^1$ by hypothesis.
Since $\|H_1(t)(H_1(0)+i)^{-1}\|$ is uniformly bounded in
$t$, we get
\bea
(H_1(t)P_0(t))'\ffi
&=&H_1(t)(H_1(0)+i)^{-1}\left((H_1(0)+i)P_0(t)\right)'\ffi\\
&+&\left(H_1(t)(H_1(0)+i)^{-1}\right)'(H_1(0)+i)P_0(t)\ffi\nonumber
\eea
where all operators are bounded and strongly continuous.\\
c)
Let $\ffi\in D_{1}$. The vector
\be
  (P_0(t)H_1(t))'\ffi={P_0}'(t)H_1(t)\ffi+P_0(t)H_1'(t)\ffi
\ee
is well defined and continuous since $P_0(t)$ is uniformly bounded in
$t$.  On the one hand, it follows from the foregoing that
\be
  {P_0}'(t)H_1(t)=(H_1(t)P_0'(t))^*
\ee
where $H_1(t)P_0'(t)$ is bounded and strongly $C^0$. Hence
${P_0}'(t)H_1(t)$ is uniformly bounded in $t$. On the other hand, as
$H_1(t)$ is self adjoint, $H_1'(t)$ is symmetric on $D_{1}$.
Thus we can write
\be
  P_0(t)H_1'(t)\subset P_0(t){H_1'}^*(t)\subset (H_1'(t)P_0(t))^*
\ee
where
\be
  H_1'(t)P_0(t)=H_1'(t)(H_1(0)+i)^{-1}(H_1(0)+i)P_0(t).
\ee
Now
\be
  H_1'(t)(H_0(0)+i)^{-1}={H_1'}^*(t)(H_1(0)+i)^{-1}
\ee
so this operator is closed and
everywhere defined, hence bounded, and it is strongly $C^0$. The same
is true for $(H_1(0)+i)P_0(t)$ so that applying the extension
principle again we eventually  obtain the uniform boundedness in $t$
of the operator
$P_0(t)H_1'(t)$. It follows that $(P_0(t)H_1(t))'$ is strongly
continuous on a dense domain and is uniformly bounded. Thus, by virtue
of theorem 3.5 p. 151 in \cite{ka2}, it is strongly continuous on
${\cal H}$.
\ep

\section{Proof of Lemma \protect\ref{GPR} b):}
\setcounter{equation}{0}

By definition
\bea
  && P_0(t){\cal R}(B)(t)P_0(t)=\nonumber\\
  &&\frac{1}{(2\pi i)^3}\oint_{\Gamma}\left(\oint_{\Gamma '}\left(
\oint_{\Gamma ''}
  R_0(t,\lambda ')R_0(t,\lambda )B(t)R_0(t,\lambda )R_0(t,\lambda '')
  d\lambda ''\right)d\lambda '\right)
  d\lambda
\eea
where the paths $\Gamma\; ,\; \Gamma '$ and $\Gamma ''$ are in the
resolvent set of $H_0(t)$, do not intersect and $\Gamma$ surrounds
$\Gamma'$ which surrounds
$\Gamma''$ which surrounds $\sigma_0^b(t)$.
We can write the integrand under the following form, using the first
resolvent equation
\bea
  \frac{R_0(t,\lambda ')B(t)R_0(t,\lambda )}{(\lambda '-\lambda )
  (\lambda -\lambda '')}-
  \frac{R_0(t,\lambda ')B(t)R_0(t,\lambda '')}{(\lambda '-\lambda )
  (\lambda -\lambda '')} \\
  -\frac{R_0(t,\lambda )B(t)R_0(t,\lambda )}{(\lambda '-\lambda )
  (\lambda -\lambda '')}+
  \frac{R_0(t,\lambda )B(t)R_0(t,\lambda '')}{(\lambda '-\lambda )
  (\lambda -\lambda '')}.\nonumber
\eea
Now, integrating each term over the variable that does not appear in
the  resolvents, we obtain the result by the Cauchy formula. For the
term $Q_0(t) {\cal R}(B)(t)Q_0(t)$ we use the definition of $Q_0(t)$
and the above result to obtain
\be
  Q_0(t){\cal R}(B)(t)Q_0(t)={\cal R}(B)(t)-P_0(t){\cal R}(B)(t)-
  {\cal R}(B)(t)P_0(t).
\ee
With the same paths as above we compute
\bea\label{pb}
  & &-P_0(t){\cal R}(B)(t)=\\
  & &\frac{1}{(2\pi i)^2}\oint_{\Gamma}\oint_{\Gamma '}
  \frac{R_0(t,\lambda ')B(t)R_0(t,\lambda )}{(\lambda '-\lambda )}-
  \frac{R_0(t,\lambda )B(t)R_0(t,\lambda )}{(\lambda '-\lambda )}
  d\lambda d\lambda '
    \nonumber
\eea
and
\bea\label{bp}
  & &-{\cal R}(B)(t)P_0(t)=\\
  & &\frac{1}{(2\pi i)^2}\oint_{\Gamma}\oint_{\Gamma '}
  \frac{R_0(t,\lambda )B(t)R_0(t,\lambda ')}{(\lambda '-\lambda )}-
  \frac{R_0(t,\lambda )B(t)R_0(t,\lambda )}{(\lambda '-\lambda )}d
\lambda
  d\lambda ' \nonumber
\eea
where the last term in (\ref{pb}) and (\ref{bp}) drops after an
integration over
$\lambda '$. Let us perform the integration over $\lambda $ in the
first term of (\ref{pb})
\bea
  & &\frac{1}{2\pi i}\oint_{\Gamma} \frac{R_0(t,\lambda ')B(t)R_0(t,
\lambda )}
  {(\lambda '-\lambda )}d\lambda= \\
  & &\frac{1}{2\pi i}\oint_{\Gamma ''} \frac{R_0(t,\lambda ')B(t)
R_0(t,\lambda '')}
  {(\lambda '-\lambda '')}d\lambda ''- R_0(t,\lambda ')B(t)R_0(t,
\lambda ')
  \nonumber
\eea
by the Cauchy formula. Thus it remains
\bea
  & &Q_0(t){\cal R}(B)(t)Q_0(t) = \\
  & &\frac{1}{(2\pi i)^2}\oint_{\Gamma '}\oint_{\Gamma ''}
\frac{R_0(t,\lambda ')B(t)
  R_0(t,\lambda '')}{(\lambda '-\lambda '')}d\lambda 'd\lambda '' -
\nonumber\\
  & &\frac{1}{(2\pi i)^2}\oint_{\Gamma}\oint_{\Gamma '}
\frac{R_0(t,\lambda )B(t)
  R_0(t,\lambda ')}{(\lambda -\lambda ')}d\lambda d\lambda '-
\nonumber \\
  & &\frac{1}{2\pi i}\oint_{\Gamma '}R_0(t,\lambda ')B(t)R_0(t,
\lambda ') d\lambda ' +{\cal R}B(t) \nonumber
\eea
where the first two terms vanish by the Cauchy formula and the last
two by definition of ${\cal R}(B)(t)$.\ep

\end{document}